\documentclass[a4paper,11pt]{article}

\usepackage{amsmath,amsthm,amssymb,mathrsfs,latexsym,amsfonts}
\usepackage{graphicx,psfrag,epsfig}
\usepackage{bbm}
\usepackage[english]{babel}
\usepackage[latin1]{inputenc}
\usepackage{a4wide}

\newtheorem{theorem}{Theorem}

\newtheorem{proposition}[theorem]{Proposition}

\newtheorem{question}{Question}

\newcounter{paraga}[section]
\renewcommand{\theparaga}{{\bf\arabic{paraga}.}}
\newcommand{\paraga}{\medskip \addtocounter{paraga}{1} 
\noindent{\theparaga\ } }

\begin{document}

\def\MP{\,{<\hspace{-.5em}\cdot}\,}
\def\SP{\,{>\hspace{-.3em}\cdot}\,}
\def\PM{\,{\cdot\hspace{-.3em}<}\,}
\def\PS{\,{\cdot\hspace{-.3em}>}\,}
\def\EP{\,{=\hspace{-.2em}\cdot}\,}
\def\PP{\,{+\hspace{-.1em}\cdot}\,}
\def\PE{\,{\cdot\hspace{-.2em}=}\,}
\def\N{\mathbb N}
\def\C{\mathbb C}
\def\Q{\mathbb Q}
\def\R{\mathbb R}
\def\T{\mathbb T}
\def\A{\mathbb A}
\def\Z{\mathbb Z}
\def\demi{\frac{1}{2}}

\begin{titlepage}
\author{Abed Bounemoura~\footnote{CNRS - CEREMADE, Université Paris Dauphine \& IMCCE, Observatoire de Paris}}
\title{\LARGE{\textbf{Non-degenerate Liouville tori are KAM stable}}}
\end{titlepage}

\maketitle

\begin{abstract}
In this short note, we prove that a quasi-periodic torus, with a non-resonant frequency (that can be Diophantine or Liouville) and which is invariant by a sufficiently regular Hamiltonian flow, is KAM stable provided it is Kolmogorov non-degenerate. When the Hamiltonian is smooth (respectively Gevrey-smooth, respectively real-analytic), the invariant tori are smooth (respectively Gevrey-smooth, respectively real-analytic). This answers a question raised in a recent work by Eliasson, Fayad and Krikorian (\cite{EFK}). We also take the opportunity to ask other questions concerning the stability of non-resonant invariant quasi-periodic tori in (analytic or smooth) Hamiltonian systems.
\end{abstract}

\section{Introduction and main result}

\paraga Let $n \geq 2$ and $\T^n:=\R^n/\Z^n$. Consider a Hamiltonian system on $\T^n \times \R^n$ associated to a $C^l$, $l \geq 2$, function of the form
\begin{equation}\label{Ham1}
H(\theta,I)=\omega \cdot I + A(\theta)I\cdot I+R(\theta,I), \quad (\theta,I)\in \T^n \times \R^n  
\end{equation}
where $\cdot$ denotes the Euclidean inner product, $\omega \in \R^n$ is a non-resonant vector ($k\cdot \omega\neq 0$ for any $k\in \Z^n \setminus\{0\}$), $A(\theta)$ is, for each $\theta \in \T^n$, a square symmetric matrix of size $n$ with real coefficients and $R(\theta,I)=O_3(I)$ is of order at least $3$ in $I$. We will also be interested in the smooth case $l=\infty$ and in the $\alpha$-Gevrey case, $\alpha \geq 1$ (with $\alpha=1$ corresponding to the real-analytic case). The set $\mathcal{T}_\omega:=\T^n \times \{I=0\}$ is invariant by the Hamiltonian flow of $H$, it is a Lagrangian quasi-periodic torus with frequency $\omega$, and any such torus (on an arbitrary symplectic manifold) is of this form. The invariant torus is said to be \emph{Kolmogorov non-degenerate} if the symmetric matrix $A_0:=\int_{\T^n}A(\theta)d\theta$ is non-singular. 

If $H$ is $C^{\infty}$, $\omega$ Diophantine (for some constant $\gamma>0$ and $\tau\geq n-1$, $|k\cdot\omega|\geq \gamma|k|_1^{-\tau}$ for any $k\in \Z^n \setminus\{0\}$, where $|k|_1:=|k_1|+\cdot+|k_n|$ if $k=(k_1,\dots,k_n)$) and $\mathcal{T}_\omega$ Kolmogorov non-degenerate, then it is \emph{KAM stable}: in any sufficiently small neighborhood of $\mathcal{T}_\omega$, there is a set of smooth Lagrangian quasi-periodic invariant tori which has positive Lebesgue measure and density one at $\mathcal{T}_\omega$. If the Hamiltonian is real-analytic, the tori are real analytic. This follows at once from a Birkhoff normal form and a classical version of the KAM theorem.

Now if $\omega$ is Liouville (which means not Diophantine), the Birkhoff normal form no longer makes sense. In \cite{EFK}, it is proved (among several other results) that if $H$ is real-analytic and $\omega$ has a ``finite uniform Diophantine exponent", the torus is still KAM stable provided it is Kolmogorov non-degenerate. The latter arithmetic condition on $\omega$, which is always satisfied for $n=2$, is satisfied for a residual subset of Liouville vectors but not for all of them if $n \geq 3$. In \cite{EFK}, the authors asked whether this arithmetic condition is necessary: in this note we prove that it is not, and that the real-analyticity assumption is also unnecessary. Informally, our main result is a follows

\begin{theorem}\label{mainthm}
The torus $\mathcal{T_\omega}$ is KAM stable provided it is Kolmogorov non-degenerate and $l$ is sufficiently large: in this case, the tori are only finitely differentiable. If $l=\infty$, the tori are smooth. If $H$ is $\alpha$-Gevrey, $\alpha\geq 1$, the tori are $\alpha$-Gevrey.
\end{theorem}

We refer to Theorem~\ref{mainthm1}, Theorem~\ref{mainthm2} and Theorem~\ref{mainthm3} below for more precise and quantitative statements concerning the regularity assumptions and conclusions, the threshold on the size of the neighborhood and the measure estimate on the set of invariant tori. Unsurprisingly, the more Liouville is the vector $\omega$, the smaller is threshold on the size of the neighborhood and the smaller is the estimate on the measure. Also, the more regular is the Hamiltonian, the better is the estimate on the measure. The proof of the main result consists, by the use of some standard scalings and some normal forms results proved in \cite{Bou13a} and \cite{Bou13b}, in reducing the above situation to a simpler situation in which classical KAM theorems apply (such as those in \cite{Pos82} or \cite{Pop04}). 

\paraga The situation becomes much more complicated for Kolmogorov degenerate invariant tori. In \cite{EFK}, there are many other interesting results concerning KAM stability of invariant tori. Perhaps the most remarkable one is that if $H$ is real-analytic and $\omega$ Diophantine, then, without any further assumptions, $\mathcal{T}_\omega$ is accumulated by KAM tori. If one assumes moreover that $n=2$ , then $\mathcal{T}_\omega$ is KAM stable (this is a result of Rüssmann) but for $n \geq 3$, it is not known if this set has positive Lebesgue measure, which was a question asked by Herman (in a related context, see \cite{Her98}). In fact, even without the Diophantine condition on $\omega$, the following question is open.

\begin{question} 
Assume $\omega$ is non-resonant and $H$ real-analytic, is $\mathcal{T}_\omega$ accumulated by KAM tori? 
\end{question}  

A positive answer to the above question would be quite surprising. However, constructing an example showing that the answer is negative seems very difficult, as even the simpler question below is open.

\begin{question} 
Assume $\omega$ is non-resonant and $H$ real-analytic, can $\mathcal{T}_\omega$ be ``unstable" in the following sense: for any $\varepsilon>0$ sufficiently small, there exists an orbit $(\theta(t),I(t))$ and a time $\tau=\tau(\varepsilon)$ such that $|I(0)|<\varepsilon$ and $|I(\tau)|>2\varepsilon$?
\end{question} 

One should expect that the answer here is positive, even for a ``generic" Hamiltonian if $n \geq 3$, but no examples are known for the moment.

If one assumes that $H$ is smooth but not analytic, an example of unstable torus is contained in \cite{Dou88}, and an example of torus which is not accumulated by a set of positive Lebesgue measure of invariant tori is contained in \cite{EFK}. The example of \cite{EFK} is, in fact, accumulated by a set invariant tori (along a hyperplane, so this set has zero measure). Hence the following question is still open.  

\begin{question} 
Assume $\omega$ is non-resonant and $H$ smooth, is $\mathcal{T}_\omega$ accumulated by KAM tori? 
\end{question}    

\section{Proof}

\paraga Let us start with the case of finitely differentiable Hamiltonians. For later use, we fix a real number $l_0>3n-1$ and we assume that $l\geq l_0+1>3n$. Consider a Hamiltonian as in~\eqref{Ham1} on the domain
\begin{equation*}
\T^n \times B_{2\varepsilon}:=\T^n \times \{I\in \R^n  \; | \; |I|<2\varepsilon \}, \quad |I|:=\max_{1 \leq i \leq n}|I_i|,
\end{equation*}
for some small $0<\varepsilon<1$. Up to the scalings
\begin{equation}\label{scale1}
I \mapsto \varepsilon I, \quad H \mapsto \varepsilon^{-1}H
\end{equation}
it is equivalent to consider the Hamiltonian
\begin{equation}\label{Ham2}
H(\theta,I)=\omega \cdot I + \varepsilon A(\theta)I\cdot I+\varepsilon^2 R(\theta,I), \quad (\theta,I)\in \T^n \times \R^n  
\end{equation}
on the domain $\T^n \times B_2$, with the estimates
\[ |A|_{C^l(\T^n)} \leq C_1, \quad |R|_{C^l(\T^n \times B_2)} \leq C_2. \]
Here, $|\,.\,|_{C^l(\T^n \times B_2)}$ denotes the usual $C^l$-norm on the domain $\T^n \times B_2$. Without loss of generality, we assume that $|\omega|=1$. For any $Q \geq 1$, let us define the function $\Psi=\Psi_\omega$ by
\[ \Psi(Q):=\max\{ |k\cdot\omega|^{-1} \; | \; k \in \Z^n, \; 0<|k|_1\leq Q \} \]
and for any $x\geq 1$, the function $\Delta=\Delta_\omega$ by
\[ \Delta(x):=\sup\{Q\geq 1 \; | \; Q\Psi(Q) \leq x\}. \]
In the sequel, we shall denote by $c_i$, for $i=1,\dots,10$, positive constants which depend only on $n$, $C_1$, $C_2$ and the operator norm of the constant matrix $A_0$ and its inverse $A_0^{-1}$, and by $c$ some universal constant which depends only on $n$.

\paraga Let us now define $f=f_\varepsilon$ by
\[ f(\theta,I):=A(\theta)I\cdot I+\varepsilon R(\theta,I) \]
so that the Hamiltonian~\eqref{Ham2} can be written as
\begin{equation*}
H(\theta,I)=\omega \cdot I + \varepsilon f(\theta,I).
\end{equation*}
Let us also define the average of the perturbation by
\[ \bar{f}(I):=\int_{\T^n}f(\theta,I)d\theta=A_0I\cdot I +\varepsilon \bar{R}(I) \]
and another ``small" parameter
\begin{equation}\label{mu}
\mu(\varepsilon):=\Delta_{\omega}^*(c\varepsilon^{-1})^{-1}.
\end{equation}
We have the following proposition.

\begin{proposition}\label{normalform}
Assume that $\mu(\varepsilon)\leq c_1$. Then there exists a symplectic embedding $\Phi : \T^n \times B_1 \rightarrow \T^n \times B_2$ of class $C^{l-1}$ such that
\[ \tilde{H}(\theta,I):=H \circ \Phi (\theta,I) = \omega \cdot I +\varepsilon \bar{f}(I)+ \varepsilon\mu(\varepsilon)\tilde{f}(\theta,I) \]
with the estimates $|\Phi-\mathrm{Id}|_{C^{l-1}(\T^n \times B_1)} \leq c_2\mu(\varepsilon)$ and $|\tilde{f}|_{C^{l-1}(\T^n \times B_1)} \leq c_3$.
\end{proposition}

This is a very special case of Theorem $1.1$ in \cite{Bou13b}, to which we refer for a proof (this proof is actually based on a result contained \cite{BF13}). To explain the difference between this (one step) normal form and (one step of) Birkhoff normal form, assume for simplicity that $f(\theta,I)$ reduces to the quadratic part $A(\theta)I\cdot I$ (the higher order terms $R(\theta,I)$ are then normalized by successive steps of Birkhoff procedure). To normalize the quadratic part, one has to solve a usual homological equation but the difficulty here (as opposed to the case of an elliptic equilibrium point) is that when expanded in Fourier series, $f(\theta,I)$ contains arbitrarily high harmonics. If $\alpha$ satisfies a Diophantine condition, the homological equation can be solved and the quadratic part normalized (assuming a Diophantine condition, this scheme would be much better as one could obtain a normal form as above but with $\varepsilon^2$ instead of $\varepsilon\mu(\varepsilon)$). When $\alpha$ is not Diophantine, the equation can not be solved exactly but only approximately. One possible approximation is to replace $f(\theta,I)$ by a trigonometric polynomial in $\theta$ (with coefficients that are functions of $I$). With the aim of proving persistence of invariant tori with a condition weaker than the Diophantine condition, Rüssmann has obtained quite precise results on the approximation of real-analytic periodic function by trigonometric polynomials (\cite{Rus01}, \cite{Rus10}) that could be applicable here (assuming real-analyticity). A different approximation scheme was proposed in \cite{BF13}: instead of approximating the perturbation $f$, the frequency $\omega$ itself is approximated by linearly independent periodic frequencies: at least one advantage of this method is that it gives quite precise results for non-analytic Hamiltonians (as explained in \cite{Bou13a} and in \cite{Bou13b}; in the real-analytic case such an approach is ``quantitatively" equivalent to the one of Rüssmann, compare \cite{BF13} and \cite{Pos11}). Let us also point out that using this scheme, the remainder term $\tilde{f}$ has no particular structure whereas if one could use Birkhoff procedure (which requires, once again, $\omega$ Diophantine) this remainder would contains terms of order at least $3$ in $I$: roughly speaking, with our procedure the graduation by degree (in $I$) is completely broke after one step of normal form, but fortunately this is not useful for the problem we are considering (as we are assuming a Kolmogorov non-degeneracy).

\paraga Now consider the Hamiltonian $\tilde{H}$ given by Proposition~\ref{normalform}, that can be written again as
\[ \tilde{H}(\theta,I)= \omega \cdot I +\varepsilon A_0I\cdot I +\varepsilon^2 \bar{R}(I)+ \varepsilon\mu(\varepsilon)\tilde{f}(\theta,I). \]
Up to the scalings
\begin{equation}\label{scale2}
H \mapsto \varepsilon^{-1}H, \quad t \mapsto \varepsilon t,
\end{equation}
it is equivalent to consider the Hamiltonian
\begin{equation}\label{Ham3}
\tilde{H}(\theta,I)= \varepsilon^{-1}\omega \cdot I + A_0I\cdot I +\varepsilon \bar{R}(I)+ \mu(\varepsilon)\tilde{f}(\theta,I).
\end{equation} 
We will apply to this Hamiltonian the classical KAM theorem for finitely differentiable Hamiltonians due to Pöschel, stated as Theorem A in \cite{Pos82}. It would be more natural to consider the term $\tilde{f}$ as the perturbation, and the other three terms of $\tilde{H}$ as the integrable part. However, Pöschel's result do require the integrable part to be real-analytic. Therefore we write
\[ \tilde{H}(\theta,I)=H_0(I)+H_1(\theta,I) \]
with
\[ H_0(I):=\varepsilon^{-1}\omega \cdot I + A_0I\cdot I, \quad H_1(\theta,I):=\varepsilon \bar{R}(I)+ \mu(\varepsilon)\tilde{f}(\theta,I). \]
The Hamiltonian $H_0$ is integrable, and since it is a polynomial it can be (trivially) extended as a real-analytic function on the domain
\[ V_\rho(B_1):=\bigcup_{I \in B_1}\{\xi \in \C^n \; | \; |\xi-I|<\rho \} \]
for any $\rho>0$ (we may simply choose $\rho=1$). Moreover, it is non-degenerate in the sense of Kolmogorov: the Hessian matrix of $H_0$ at any point is given by $2A_0$ and so it is non-singular. Since $\varepsilon$ is (much) smaller than $\mu(\varepsilon)$, the size of the perturbation $H_1$ is (up to a constant) bounded by $\mu(\varepsilon)$, and as $H_1$ is of class $C^{l-1}$ with $l-1\geq l_0>3n-1$, an application of \cite{Pos82} yields the following proposition.

\begin{proposition}\label{KAM}
Assume that $\sqrt{\mu(\varepsilon)} \leq c_4$. Then there exists a set $\tilde{\mathcal{K}}=\tilde{\mathcal{K}}_\varepsilon$ contained in $\T^n \times B_1$ which consists of Lagrangian quasi-periodic tori, invariant by the Hamiltonian flow of $\tilde{H}$. Each invariant torus is of class $C^{l_0'}$, for $l_0'<l_0-2n+1$, and we have the Lebesgue measure estimate
\[ c_5\sqrt{\mu(\varepsilon)}\mathrm{Leb}(\T^n \times B_1) \leq \mathrm{Leb}(\T^n \times B_1 \setminus \tilde{\mathcal{K}}) \leq  c_6\sqrt{\mu(\varepsilon)}\mathrm{Leb}(\T^n \times B_1). \]
\end{proposition}   

Let us make some comments on the above proposition. The smoothness requirement and conclusion are apparently different than those of \cite{Pos82}; however, as pointed out by Sevryuk in \cite{Sev03}, they are in fact the same. More importantly, since $\mu(\varepsilon)$ is much smaller than $\varepsilon$, one should explain why the $\varepsilon$-dependence of the integrable part does not obstruct the result. Observe that the frequency domain, that is the image of $B_1$ by the gradient of $H_0$, is simply the image of $B_1$ by the linear map $2A_0$ translated by the vector $\varepsilon^{-1}\omega$. So the frequencies are high (they have large norm) but the important point is that the Hessian matrix of $H_0$, which is nothing but the matrix $2A_0$, and its inverse, have bounds uniform with respect to $\varepsilon$. The set of frequencies preserved are, as usual, those which are $(\gamma,\tau)$-Diophantine, with some $\gamma>0$ and $\tau>n-1$, and which are at a distance at least $\gamma$ from the boundary; but recalling the scaling in time this actually corresponds to the persistence of (unperturbed) quasi-periodic motions with frequencies $\alpha+2\varepsilon A_0I$ which are $(\varepsilon\gamma,\tau)$-Diophantine. Eventually, the threshold is, as usual, of the form $\mu(\varepsilon)\leq c_4 \gamma^2$ for $\gamma\leq 1$, which allows to choose $\gamma$ proportional to $\sqrt{\mu(\varepsilon)}$ to obtain the measure estimate. Let us also note, for later use, that by construction, this Kolmogorov set $\tilde{\mathcal{K}}$ is at distance at least, up to a constant, $\sqrt{\mu(\varepsilon)}$ to the boundary of $\T^n \times B_1$.

\paraga Undoing the scalings~\eqref{scale2}, one finds a set of KAM tori $\tilde{\mathcal{K}}$, with the properties stated in Proposition~\ref{KAM}, which is invariant by the Hamiltonian flow of $\tilde{H}=H\circ \Phi$ given by Proposition~\ref{normalform}. Using the fact that $\mu(\varepsilon)$ is smaller than $\sqrt{\mu(\varepsilon)}$, the threshold and the estimate on $\Phi$ given by Proposition~\ref{normalform}, one easily ensures that the image of $\Phi$ contains $\tilde{\mathcal{K}}$. Therefore we can define $\mathcal{K}=\Phi^{-1}(\tilde{\mathcal{K}})$, and since $\Phi$ is symplectic and of class $C^{l_0}$, $\mathcal{K}$ is a set which consists of Lagrangian quasi-periodic tori, invariant by the Hamiltonian flow of $H$ defined by~\eqref{Ham2}, and which are of class $l_0'$, for $l_0'<l_0-2n+1$. For the measure estimate, observe that the estimate on $\Phi$ given by Proposition~\ref{normalform} implies that its Jacobian is close to one, hence 
\[ c_7\sqrt{\mu(\varepsilon)}\mathrm{Leb}(\T^n \times B_1) \leq \mathrm{Leb}(\T^n \times B_1 \setminus \mathcal{K}) \leq  c_8\sqrt{\mu(\varepsilon)}\mathrm{Leb}(\T^n \times B_1). \]
Undoing the scalings~\eqref{scale1}, one finds a set, that we sill denote by $\mathcal{K}$ for simplicity, which is contained in $\T^n \times B_{\varepsilon}$ and invariant by the Hamiltonian flow of $H$ defined by~\eqref{Ham1}, and which consists of Lagrangian quasi-periodic invariant tori, with the measure estimate 
\[ c_9\sqrt{\mu(\varepsilon)}\mathrm{Leb}(\T^n \times B_{\varepsilon}) \leq \mathrm{Leb}(\T^n \times B_{\varepsilon} \setminus \mathcal{K}) \leq  c_{10}\sqrt{\mu(\varepsilon)}\mathrm{Leb}(\T^n \times B_{\varepsilon}). \]
We have just proved the following statement.

\begin{theorem}\label{mainthm1}
Let $H$ be as in~\eqref{Ham1}, assume it is of class $C^l$ for $l\geq l_0+1>3n$, that $\omega$ is non-resonant and $A_0$ is non-singular. Then if $\mu(\varepsilon) \leq c_1$ and $\sqrt{\mu(\varepsilon)} \leq c_4$, where $\mu(\varepsilon)$ is defined in~\eqref{mu}, there exists a set $\mathcal{K} \subset \T^n \times B_\varepsilon$, which consists of Lagrangian quasi-periodic tori invariant by the Hamiltonian flow of $H$. Moreover, each tori is of class $C^{l_0'}$, for $l_0'<l_0-2n+1$, and we have the measure estimate
\[ c_9\sqrt{\mu(\varepsilon)}\mathrm{Leb}(\T^n \times B_\varepsilon) \leq \mathrm{Leb}(\T^n \times B_\varepsilon \setminus \mathcal{K}) \leq  c_{10}\sqrt{\mu(\varepsilon)}\mathrm{Leb}(\T^n \times B_\varepsilon). \]
\end{theorem}

This result justifies the first part of Theorem~\ref{mainthm}. As a matter of fact, since Theorem~\ref{mainthm1} reduces to a classical KAM theorem, more information is available. For instance, this set $\mathcal{K}$, which forms a Cantor family, is in fact regular (in the sense of Whitney) over this Cantor set: here this transverse regularity is $l_0''$, for any $l_0''<(l_0-2n+1)/n$. But of course only a Lipschitz transverse regularity is necessary to obtain the measure estimate. 

\paraga  When $H$ is $C^{\infty}$, using the smooth version of Theorem 1.1 in \cite{Bou13b} and the smooth version of Theorem A in \cite{Pos82}, one gets the following result.

\begin{theorem}\label{mainthm2}
Assume that $H$ is $C^{\infty}$, and that the assumptions of Theorem~\ref{mainthm1} are satisfied. Then the conclusions of Theorem~\ref{mainthm1} hold true, and, in addition, the tori are $C^{\infty}$.
\end{theorem}

This justifies the second part of Theorem~\ref{mainthm}. Here one would expect more conclusions: for instance, given any fixed integer $\kappa\geq 1$, the measure of the set not covered by KAM tori should be of order $\sqrt{\mu(\varepsilon)^\kappa}$. Indeed, assuming $H$ is smooth, Theorem 1.1 in \cite{Bou13b} gives a normal form 
\[ \tilde{H}(\theta,I):=H \circ \Phi (\theta,I) = \omega \cdot I +\varepsilon \bar{f}(I)+\varepsilon\mu(\varepsilon)g(I)+ \varepsilon\mu(\varepsilon)^\kappa\tilde{f}(\theta,I) \]
where $g$ is integrable (for $\kappa=1$, which is the case considered in Proposition~\ref{normalform}, this Hamiltonian $g$ can be taken to be identically zero). After scaling as in~\eqref{scale2}, the dominant part in the integrable part of $\tilde{H}$ is still $\bar{f}$, hence this integrable Hamiltonian is Kolmogorov non-degenerate. However, Pöschel's result \cite{Pos82} demands a real-analytic integrable part: in our case it is only smooth and since $g$ is quite arbitrary, we do not know how to write $\tilde{H}$ as a real-analytic integrable Hamiltonian plus a perturbation of order $\mu(\varepsilon)^\kappa$. Let us remark here that the real-analyticity assumption in Pöschel's result seems quite artificial, so in the end one should have this improved measure estimate but of course we do not claim such a result.

\paraga To conclude, let us consider the $\alpha$-Gevrey case, for $\alpha\geq 1$, which includes the real-analytic case $\alpha=1$. More precisely, we assume that $H$ is $(\alpha,L)$-Gevrey for some $L>0$ in the sense of~\cite{MS02}. When $\alpha=1$, $H$ is real-analytic and this constant $L$ is comparable to a width of analyticity. 

Let us denote by $\bar{c}$ some universal constant depending only on $n$, $\alpha$ and $L$, and let 
\begin{equation}\label{nu}
\nu(\varepsilon):=\exp\left(-\bar{c}\mu(\varepsilon)^{-\alpha^{-1}}\right).
\end{equation}  
Using Theorem 1.1 of \cite{Bou13a} (which uses results from \cite{BF13} and \cite{MS02}), for $\mu(\varepsilon)$ small, one finds a normal form 
\[ \tilde{H}(\theta,I):=H \circ \Phi (\theta,I) = \omega \cdot I +\varepsilon \bar{f}(I)+\varepsilon\mu(\varepsilon)g(I)+ \varepsilon\nu(\varepsilon)\tilde{f}(\theta,I) \]
such that $g$ and $\tilde{f}$ are bounded in $(\alpha,\tilde{L})$-Gevrey norm, for some fixed $\tilde{L}<L$. After scalings as in~\eqref{scale2}, we do have a perturbation of order $\nu(\varepsilon)$ of a Kolmogorov non-degenerate integrable Hamiltonian: both the integrable Hamiltonian and the perturbation are $(\alpha,\tilde{L})$-Gevrey, the Hessian of the integrable part at any point is close to $2A_0$ and this is enough to apply the main theorem of \cite{Pop04}, which is an extension of Pöschel's result in Gevrey classes. Therefore we can state the following result.

\begin{theorem}\label{mainthm3}
Let $H$ be as in~\eqref{Ham1}, assume it is $(\alpha,L)$-Gevrey, that $\omega$ is non-resonant and $A_0$ is non-singular. Then if $\mu(\varepsilon) \leq c_{11}$ and $\sqrt{\nu(\varepsilon)} \leq c_{12}$, where $\mu(\varepsilon)$ is defined in~\eqref{mu} and $\nu(\varepsilon)$ is defined in~\eqref{nu}, there exists a set $\mathcal{K} \subset \T^n \times B_\varepsilon$, which consists of Lagrangian quasi-periodic tori invariant by the Hamiltonian flow of $H$. Moreover, each tori is of class $(\alpha,\tilde{L}')$, for some $\tilde{L}'<\tilde{L}$, and we have the measure estimate
\[ c_{13}\sqrt{\nu(\varepsilon)}\mathrm{Leb}(\T^n \times B_\varepsilon) \leq \mathrm{Leb}(\T^n \times B_\varepsilon \setminus \mathcal{K}) \leq  c_{14}\sqrt{\nu(\varepsilon)}\mathrm{Leb}(\T^n \times B_\varepsilon). \]
\end{theorem}

This justifies the last part of Theorem~\ref{mainthm}. Let us remark here that the $(\alpha,L)$-Gevrey norm used by \cite{Pop04} is different that the one used in \cite{MS02}; however, up to changing $L$ by its inverse, they are comparable. The above positive constants $c_{11},c_{12},c_{13},c_{14}$ depend, as before, on $n$, $C_1$, $C_2$ (which are bounds, respectively, on the $(\alpha,L)$-Gevrey norms of $A$ and $R$ that appears in~\eqref{Ham2}), the operator norm of the constant matrix $A_0$ and its inverse $A_0^{-1}$, but also on $\alpha$ and $L$. The proof of Theorem~\ref{mainthm3} goes exactly as the proof of Theorem~\ref{mainthm1}, with only one mild difference: $\sqrt{\nu(\varepsilon)}$ plays the role of $\sqrt{\mu(\varepsilon)}$, but this time it is much smaller than $\mu(\varepsilon)$. Therefore in the Gevrey version of Proposition~\ref{KAM} (which follows from~\cite{Pop04}), one looks for a set $\tilde{\mathcal{K}}$ of KAM tori not in $\T^n \times B_1$ but on smaller domain, say in $\T^n \times B_{1/2}$, so that one can easily ensures that this set is contained in the image of the transformation $\Phi$ given by the Gevrey version of Proposition~\ref{normalform} (which follows from~\cite{Bou13a}), which in turns ensures that the set of KAM tori $\mathcal{K}=\Phi^{-1}(\tilde{\mathcal{K}})$ is well-defined. This only affects the measure estimate by constants depending on $n$.    

In the special case where $\alpha$ is $(\gamma,\tau)$ Diophantine, then $\mu(\varepsilon)$ is of order $\varepsilon^{(1+\tau)^{-1}}$ and hence $\nu(\varepsilon)$ is order $\exp\left(-\varepsilon^{-(\alpha(1+\tau))^{-1}}\right)$. In the analytic case $\alpha=1$, it is well-known that the measure of the complement of the Kolmogorov set is exponentially small; our result extends this for any $\alpha \geq 1$. 

\addcontentsline{toc}{section}{References}
\bibliographystyle{amsalpha}
\bibliography{LiouvilleTori}

\providecommand{\bysame}{\leavevmode\hbox to3em{\hrulefill}\thinspace}
\providecommand{\MR}{\relax\ifhmode\unskip\space\fi MR }
\providecommand{\MRhref}[2]{%
  \href{http://www.ams.org/mathscinet-getitem?mr=#1}{#2}
}
\providecommand{\href}[2]{#2}
\begin{thebibliography}{Bou13b}

\bibitem[BF13]{BF13}
A.~Bounemoura and S.~Fischler, \emph{A diophantine duality applied to the {KAM}
  and {N}ekhoroshev theorems}, Math. Z. \textbf{275} (2013), no.~3, 1135--1167.

\bibitem[Bou13a]{Bou13a}
A.~Bounemoura, \emph{Normal forms, stability and splitting of invariant
  manifolds {I}. {G}evrey {H}amiltonians}, Regul. Chaotic Dyn. \textbf{18}
  (2013), no.~3, 237--260.

\bibitem[Bou13b]{Bou13b}
\bysame, \emph{Normal forms, stability and splitting of invariant manifolds
  {II}. {F}initely differentiable {H}amiltonians}, Regul. Chaotic Dyn.
  \textbf{18} (2013), no.~3, 261--276.

\bibitem[Dou88]{Dou88}
R.~Douady, \emph{Stabilité ou instabilité des points fixes elliptiques}, Ann.
  Sci. Ec. Norm. Sup. \textbf{21} (1988), no.~1, 1--46.

\bibitem[EFK]{EFK}
L.H. Eliasson, B.~Fayad, and R.~Krikorian, \emph{Around the stabiity of {KAM}
  tori}, Duke Mathematical Journal, To appear.

\bibitem[Her98]{Her98}
M.~Herman, \emph{Some open problems in dynamical systems}, Doc. Math., J. DMV,
  Extra Vol. ICM Berlin 1998, vol. {II}, 1998, pp.~797--808.

\bibitem[MS02]{MS02}
J.-P. Marco and D.~Sauzin, \emph{Stability and instability for {G}evrey
  quasi-convex near-integrable {H}amiltonian systems}, Publ. Math. Inst. Hautes
  Études Sci. \textbf{96} (2002), 199--275.

\bibitem[Pop04]{Pop04}
G.~Popov, \emph{K{AM} theorem for {G}evrey {H}amiltonians}, Erg. Th. Dyn. Sys.
  \textbf{24} (2004), no.~5, 1753--1786.

\bibitem[Pös82]{Pos82}
J.~Pöschel, \emph{Integrability of {H}amiltonian systems on {C}antor sets},
  Comm. Pure Appl. Math. \textbf{35} (1982), no.~5, 653--696.

\bibitem[Pös11]{Pos11}
\bysame, \emph{K{A}{M} à la {R}}, Regul. Chaotic Dyn. \textbf{16} (2011),
  no.~1-2, 17--23.

\bibitem[Rüs01]{Rus01}
H.~Rüssmann, \emph{Invariant tori in non-degenerate nearly integrable
  {H}amiltonian systems}, Regul. Chaotic Dyn. \textbf{6} (2001), no.~2,
  119--204.

\bibitem[Rüs10]{Rus10}
\bysame, \emph{{KAM}-iteration with nearly infinitely small steps in dynamical
  systems of polynomial character}, Discrete Contin. Dyn. Syst. Ser. S
  \textbf{3} (2010), no.~4, 683--718.

\bibitem[Sev03]{Sev03}
M.~B. Sevryuk, \emph{{The classical KAM theory at the dawn of the Twenty-First
  Century}}, Mosc. Math. J. \textbf{3} (2003), no.~3, 1113--1144.

\end{thebibliography}

\end{document}